\documentclass[11pt]{article}
 \usepackage{amsmath,amssymb}
 \usepackage{epsfig}
 \usepackage{graphicx}
\usepackage{verbatim}

 \newcommand{\sfrac}[2]{{\textstyle\frac{#1}{#2}}}

\newcommand{\cost}{\mathbf{cost}}

\def\Ex{\mathbb{E}}
\renewcommand{\Pr}{\mathbb{P}}

\newcommand{\var}{\mathrm{var\ }}

 \newcommand{\wins}{\mathrm{ wins}}
 
\begin{document}
 
\title{Gambling under unknown probabilities as a proxy for real world decisions under uncertainty.}
 \author{David J. Aldous\thanks{Department of Statistics,
 367 Evans Hall \#\  3860,
 U.C. Berkeley CA 94720;  aldous@stat.berkeley.edu;
  www.stat.berkeley.edu/users/aldous.}
\and F. Thomas Bruss\thanks{Facult\'{e} des Sciences, Unversit\'{e} Libre de Bruxelles, Campus Plaine, CP 210, B-1050 Brussels, Belgium; tbruss@ulb.ac.be.}
}

\date{December 8, 2020}
 \maketitle
 
 {\bf This was a draft.  The final version is published as {\em Amer. Math. Monthly} 130 (2023) 303--320.}
 
 \begin{abstract}
 We give elementary examples within a framework for studying decisions under uncertainty where probabilities
 are only roughly known.  The framework, in gambling terms, is that the size of a bet is proportional to the gambler's 
 perceived advantage based on their perceived probability, and their accuracy in estimating true probabilities is measured by mean-squared-error.  Within this framework one can study the cost of estimation errors, 
 and seek to formalize the ``obvious" notion that in
 competitive  interactions between agents whose actions depend on their perceived  probabilities, 
those who are more accurate at estimating probabilities will generally be more successful than those 
who are less accurate.
%{\tt yyy
% I guess these ideas must have been thought about before, but I haven't found any similar elementary articles.}
 \end{abstract}

 \section{Introduction}

The general topic of decisions under uncertainty covers a broad spectrum, 
from the classical mathematical decision theory surrounding expected utility \cite{decision1}
to modern work (with the celebrated popular exposition by Kahneman \cite{kahneman})
on the cognitive biases exhibited by  actual human beings in making such decisions.
An important practical point for any discussion is that, for perhaps the majority of real-world decisions
we face outside specific professional contexts, numerical probabilities
can only be guessed or crudely estimated.
Where theory and data do not enable us to estimate probabilities well, can one still seek to 
study the accuracy of human-estimated probabilities and the effects of errors in such estimation?

Such questions have been discussed at sophisticated levels in various specific contexts 
(see section \ref{sec:sophist} discussion) but this article arises from thinking how to introduce the topic 
-- what could one say in a single class as part of an undergraduate course in probability?

To pose a starting question, consider the following two intuitive notions.
\begin{quote}
(A) For unique real-world events, such as ``will this politician be re-elected next year?", either there is no ``true probability", or there is some unknown true probability and we can never tell in any quantitative sense
whether your guess of a 60\% chance was better than my guess of 40\% chance.
\end{quote}
On the other hand, in sports betting or stock market speculation, some individuals do better than others, 
perhaps more than by pure chance, and one can formulate
%In an attempt to think more generally about such questions, consider 
the following vague general assertion.
\begin{quote}
(B) In any kind of competitive  interactions ``under uncertainty" between agents whose actions depend on their perceived  probabilities, 
those who are more accurate at estimating probabilities will be more successful than those 
who are less accurate.
\end{quote}
Notions (A) and (B) are not exactly contradictory, but represent ends of a spectrum of views about 
assessment of numerical probabilities for interesting future real-world events. 

To the extent that  the effects of different outcomes can be expressed within the same quantitative units 
 (as assumed in utility theory), any decision under uncertainty model can be regarded as a gambling model.
This article outlines a framework 
within which assertion (B) above, 
interpreted as involving ``gambling under unknown probabilities", can be studied mathematically.
We could start with a detailed discussion of background and motivation for this framework.
But anticipating that mathematical readers are likely to be frustrated by many pages of verbal discussion without many $x$'s and $y$'s, let us first just outline the framework and then
jump quickly into a study of 6 toy models, whose analysis involves only straightforward mathematics.
So further general discussion 
and pointers to the academic literature are deferred to  sections \ref{sec:discuss} and \ref{sec:sophist}.

\subsection{Outline of article}
\begin{figure}
\caption{Stories}
\label{Fig:pics}

\includegraphics[width=1.1in]{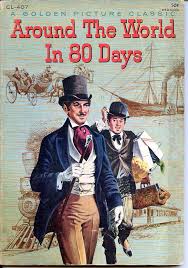} \hspace*{0.3in} \vspace*{0.1in}
\includegraphics[width=1.4in]{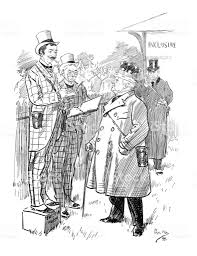}  \hspace*{0.3in} 
\includegraphics[width=1.6in]{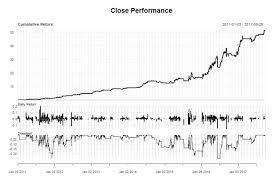}

\includegraphics[width=1.5in]{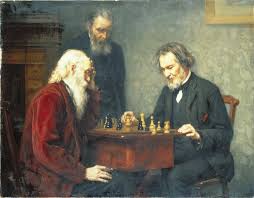}
\includegraphics[width=1.6in]{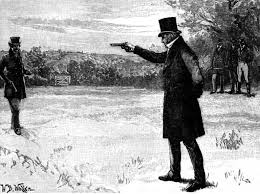}
\includegraphics[width=1.6in]{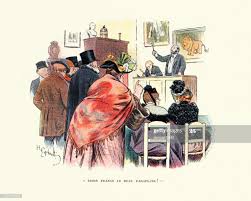}
\end{figure}

%We will explain the background and setup slowly in sections \ref{sec:cost} and \ref{sec:frame}, but here is our framework in brief.
Here is the basic framework.
\begin{quote}
Take a toy model of a situation where one has to make an action (like deciding whether and how much to bet)
whose outcome (gain/loss of money/utility) depends on whether an event of probability $p_{true}$ occurs.
There is some known optimal (maximize expected utility) action if $p_{true}$ is known.
But all one has is a ``perceived" 
probability $p_{perc}$.  So one just takes the action that one would take if $p_{perc}$ were the true probability.
Now we study the consequences of the action under the assumption 
 that $p_{perc} = p_{true} + \xi$ for random error $\xi$,
where usually we need to assume that $\xi$ has mean zero.
\end{quote}
Subsequent sections give analyses of different models, indicated informally here and illustrated in Figure \ref{Fig:pics}.
\begin{itemize}
\item  {\bf A gentleman's bet}: You think a future event has probability 20\%, your friend thinks it has probability 30\%, so you make a bet at odds corresponding to 25\%. (section \ref{sec:2pers})
\item {\bf  The bookmakers dilemma}:  A bookmaker offers odds corresponding to different event probabilities, 
say $64\%$ and $60\%$, for an event happening or not happening.  How to choose these values, based on the bookmakers and the gamblers' perceptions of the probability? (section \ref{sec:book})
\item {\bf Bet I'm better than you!}: 
Two opponents in a game may choose to bet at even odds, but only do so if each believes they are more skillful than the other.
(section \ref{sec:skill})
\item {\bf Kelly rules}: Adapting to our setting the Kelly criterion for allocating sizes of favorable bets. 
(section \ref{sec:kelly})
\end{itemize}
The models above fit into the
basic setting where the only unknown quantity is the probability of a given event.
The following models have more elaborate settings of ``unknowns".
\begin{itemize}
\item {\bf Pistols at dawn}: When to fire your one shot, if uncertain about abilities. (section \ref{sec:pistols})
\item {\bf How valuable is it really?}:  Unknown utilities when choosing or bidding.  (section  \ref{sec:Gumbel})
\end{itemize}
The attentive reader may wonder about our default assumption that agents act as if their perceived probability were the true probability -- could they instead
try to make allowance for likely error?
This {\em allowance issue} is a fascinating question for future study but we will only address it  occasionally in this article.

\subsection{Other model ingredients}
There are 3 more ingredients in almost all our examples.
First, in gambling-like settings, it is convenient to 
use the terminology (from prediction markets)
% discussed in more detail in section xxx)
 of {\em contracts} rather than {\em odds}.
A contract on an event will pay \$1 if the event occurs, zero if not.
In traditional horse race language one might bet \$7 at odds of 4-to-1 against;
that means you would  gain \$28 if the horse wins, or lose \$7 if not.
This corresponds, in our terminology, to buying 35 contracts at price $0.20$ dollars per contract.
Here the price is the {\em implied probability} $p_{implied} = 0.20$ if these were ``fair odds". 

Second, we need to model the amount that is bet in any particular case.
Of course in most circumstances
% it is basically rational (yyy discuss at end) 
you would bet on that horse only if 
 your perceived\footnote{Terminology is awkward: {\em estimated} suggests an explicit estimation rule, {\em subjective} suggests it's just an opinion.  So {\em perceived} is a compromise.}  probability  $p_{perc}$ 
 is greater than $p_{implied}$.
Intuition suggests that you should bet more if your perceived advantage $p_{perc} - p_{implied}$ is larger,
and we  will use the  simplest implementation of that intuition by modelling that the 
number of contracts bought is proportional to  $p_{perc} - p_{implied}$. 
This is roughly the {\em Kelly strategy} (section  \ref{sec:kelly}).
 So the number of contracts bought is $\kappa (p_{perc} - p_{implied})$, where the ``constant of proportionality" 
 $\kappa$, measuring the scale of an individual's gambling budget (``affluence", say),
   is unimportant for our analyses.
 
 Third, recall our ``unknown true probability" $p_{true}$ framework.
 Any instance of a bet by an individual involves a triple like $(p_{implied}, p_{perc}, p_{true})$ above,
and then our model 
\[ \mbox{ $p_{perc} = p_{true} + \xi$ for random error $\xi$} \]
will allow us to study mean outcomes in terms of the distribution of the perception error $\xi$.

This completes a framework for studying the ``obvious" vague assertion (B):
is it true that, other things being equal, agents with smaller error $\xi$ have better outcomes? 
Often we will need to make the ``unbiased" assumption that the error $\xi$ is such that $\Ex \xi = 0$ -- 
see section \ref{sec:calibrate} for discussion.
Writing $\sigma^2 = \var \xi$, the ``accuracy of $p_{perc}$" can be expressed as the RMS 
(root mean square) error $\sigma$.

\subsection{The error-squared principle}
It turns out, perhaps unsurprisingly, that in many cases the quantitative effect of inaccuracy of perceived probabilities 
scales as $\sigma^2$, and we will observe this in our toy model analyses.
In fact this general {\em error-squared principle}
is well known in the related context of {\em prediction tournaments} discussed in section \ref{sec:rPT},
one of the few contexts where we have direct data on human ability at estimating probabilities.

Finally, throughout this article we are {\bf not} envisaging unlikely events with large consequences, for which 
squared error is clearly not the appropriate measure.
In other words, we are implicitly assuming that $p_{true}$ is not close to $0$ or $1$.
This assumption was also used in identifying  ``bet proportional to perceived advantage" as approximately the Kelly strategy.

\section{A gentleman's bet}
\label{sec:2pers}
For our first toy model, let us check that our framework gives a plausible answer in the simplest setting.
Suppose two people (A and B) have different perceived probabilities $q_A$ and $q_B$
for a future event and wish to make a bet.
Then a contract (to receive 1 dollar if the event occurs) at any price between $q_A$ and $q_B$
is perceived as favorable to each person.
Suppose $q_A > q_B$, so A will buy and B will sell, and suppose the price is set\footnote{{\em Gentleman}, or {\em Lady},  has a certain significance here; if one
were not a gentleman then one might consider whether mis-stating your perceived probability would be advantageous.}
at the 
midpoint $(q_A + q_B)/2 := r$.
In terms of the unknown true probability $p$, the mean gain to $A$ per contract is $p - r$,
so the default strategy for A is to buy $\kappa (q_A - r)$ contracts, and then
\begin{equation}
\mbox{ mean gain to A = $\kappa (q_A - r)(p-r)$.}
\label{format}
\end{equation}
This formula also holds in the case $q_A < q_B$, where A sells $\kappa (r- q_A )$ contracts 
with  a  mean gain  per contract of $r - p$.

Within our model $q_A$ and $q_B$ are random, so
\[
\Ex [ \mbox{ mean gain to A }] = \kappa \Ex [(q_A - \sfrac{q_A +  q_B}{2}) (p -  \sfrac{q_A +  q_B}{2} ) ].
\]
The {\em unbiased} assumption is that $\Ex q_A = \Ex q_B = p$.
Writing $\sigma_A^2$ and $\sigma_B^2$ for the variances of 
$q_A$ and $q_B$, we find
\begin{equation}
\Ex [ \mbox{ mean gain to A }] = \sfrac{1}{4} \kappa (\sigma^2_B - \sigma_A^2) 
\label{eq:2-person}
\end{equation}
{\em even if $q_A$ and $q_B$ are dependent}.

So, at least here, our framework does give explicit and intuitive results:
the ``error-squared principle" is visible in a very simple way.
Note the format we used above:  equation (\ref{format}) refers to the randomness of event outcome, in terms 
of perceived and true probabilities; then we use some model of how the perceived and true probabilities are related to calculate the expectation (\ref{eq:2-person})  under the true probability within that model.  
We will use the same format in other examples.

\section{The bookmakers dilemma}
\label{sec:book}
\begin{quote}
Bookmakers are more skilled at predicting the outcomes of games than bettors and systematically exploit bettor biases by choosing prices that deviate from the market clearing price.
{\em Levitt} \cite{levitt}.
\end{quote}

\begin{quote}
Sports betting is the only [casino] game where you are, in fact, playing against the house ....... 
You're playing against other people who are actively trying to beat you.
{\em Miller - Davidow} \cite{sports_betting}
\end{quote}
Levitt  \cite{levitt} gives a detailed account of actual bookmaker strategy in the U.S. sports gambling context,
and Miller - Davidow \cite{sports_betting} describes the nuts and bolts of sports betting from the gambler's viewpoint.
Our model below is very crude -- let's see whether it is qualitatively reasonable. 

In our setup an idealized bookmaker announces a bid price $x_1$ for a gambler wishing to sell a contract, and an ask price $x_2 > x_1$ for a  gambler wishing to buy.
That is, in the context on betting whether team $T$ will win an upcoming game, 
you can bet on ``win" by paying $x_2$ and receiving $1$ if $T$ wins, or you can bet on ``lose" by 
paying $1 - x_1$ and receiving $1$ if $T$ loses.\footnote{This is equivalent to paying $1 - x_1$ for a contract for the opponent to win.}
In our model a gambler with perceived probability $p_{perc} > x_2$ will buy $\kappa (p_{perc} - x_2)$ contracts,
and a gambler with perceived probability $p_{perc} < x_1$ will sell $\kappa (x_1 - p_{perc})$ contracts.
How should the bookmaker choose the {\em spread  interval} $[x_1,x_2]$?
If the interval is too wide,\footnote{In U.S. sports betting the width $x_2-x_1$ of the spread interval is called the {\em hold} or the {\em vig} or the {\em juice} or the {\em take} or the {\em house cut}
 \cite{sports_betting}.} 
 fewer gamblers will bet, whereas if the interval is too narrow then the bookmaker 
may make too little average profit per bet.
Note that we are implicitly assuming a ``monopoly" or cartel of bookmakers, so that gamblers have no alternate venue; otherwise there would be other factors arising from competition between bookmakers.\footnote{Gamblers of course want a small spread interval.  Real-world bookmakers present odds in a variety of ways, and 
(as \cite{sports_betting} emphasizes) a prospective gambler should first learn to translate into our {\em implied odds} format and calculate the spread.}

Assume that for a given event the bookmaker knows the distribution (over gamblers) of perceived probabilities,
and that this distribution\footnote{One could analyze other distributions.} is uniform over some interval, 
so we can write this interval as 
$[p_{gamb} - L,p_{gamb}  + L]$, so that $p_{gamb}$ is the consensus probability amongst gamblers.
So for a given event\footnote{$\kappa$ here is the sum of the individual gamblers' affluences $\kappa$,
assumed independent of their perceived probabilities.}
\begin{eqnarray}
\lefteqn{\mbox{ mean gain to bookmaker   }} &&\nonumber \\
&=&  \kappa (x_2 - p_{true})   \sfrac{1}{2L}   \int_{x_2}^{p_{gamb}  + L} (x-x_2) dx 
+  \kappa (p_{true}-x_1)   \sfrac{1}{2L}   \int_{p_{gamb}  -L}^{x_1}  (x_1 -x) dx \nonumber \\
&=& \sfrac{\kappa}{4L} 
\left[  (x_2 - p_{true}) (p_{gamb}  + L - x_2)^2 + (p_{true}-x_1)  (x_1 - p_{gamb} + L)^2
\right] .
\label{Mbook1}
\end{eqnarray}
Note this is assuming 
\begin{equation}
 [x_1,x_2] \subseteq [p_{gamb} - L,p_{gamb}  + L];
 \label{consist1}
 \end{equation}
 which is reasonable because a bookmaker obviously prefers gamblers to have a wide range of  perceived probabilities to encourage actual betting, as will be seen in (\ref{Mbook2}, \ref{Mbook3})  below.  

{\bf In the case where the bookmaker knows $p_{true}$,} the bookmaker can
optimize (\ref{Mbook1}) over $x_1$ and $x_2$, and the optimal spread interval is
\begin{equation}
   [x_1,x_2] = [\sfrac{2}{3} p_{true}  + \sfrac{1}{3} (p_{gamb}  - L)      , \sfrac{2}{3} p_{true} + \sfrac{1}{3} (p_{gamb}  + L) ] 
 \label{Mbook0}
 \end{equation}
and the resulting profit is
\begin{eqnarray}
\Ex  [ \mbox{ mean gain to bookmaker   (known $p_{true}$) } ] &=&
\sfrac{2 \kappa}{27} (L^2 + 3 \Delta^2); \label{Mbook2} \\
 \Delta &:=& p_{gamb} - p_{true} . \nonumber
\end{eqnarray}
This is assuming   $p_{true} \in [p_{gamb} - L,p_{gamb}  + L]$, to satisfy the constraint (\ref{consist1}).
 
{\bf In an opposite case where the bookmaker does not know or wish to estimate $p_{true}$}, the bookmaker  could 
instead take the spread interval to be a symmetric interval around $p_{gamb}$.
Applying (\ref{Mbook1}) with $[x_1,x_2] = [p_{gamb} - x, p_{gamb} + x]$
(or calculating directly)
we find that the gain is maximized at $x = \sfrac{1}{3}L$ and the maximized value is
\begin{equation}
\Ex  [ \mbox{ mean gain to bookmaker   (interval centered at $p_{gamb}$) } ] =
\sfrac{2 \kappa}{27} \ L^2 . \label{Mbook3}
\end{equation}
Results (\ref{Mbook2}, \ref{Mbook3}) are qualitatively what one would expect.
Order $L^2$ arises as the variance of the gamblers' perceived probabilities, 
and order $\Delta^2$ as the effect of the gambler's bias when the bookmaker is not biased. 
Note the error-squared principle again.
Note also that the interval in the first case is not symmetric about either
$p_{true}$ or $p_{gamb}$, but rather about a weighted average.

\subsection{Using the bookmaker's perceived probability.}
To assume (as above)  that the bookmaker knows the true probability is unrealistic.
Instead, can we
study a model where, as in the Levitt quote above, the bookmaker is only {\em more accurate} than 
the gamblers at assessing probabilities?
In our framework we model the bookmaker's perceived probability as
\[ p_{book} =  p_{true} + \xi \]
where here the error $\xi$ has a symmetric distribution with variance $\sigma^2$.
This leads to our first engagement with the {\em allowance principle}.
Let us continue assuming that the gamblers act as if their perceived probabilities were true probabilities, but let us not require that the bookmaker does so.
Note that the bookmaker could do so, that is
could use the interval (\ref{Mbook0}) with $p_{book}$ in place of $p_{true}$.
But is this optimal in this model?

In this section we will calculate the optimal spread interval $[x_1,x_2]$ in terms of  $(p_{book},L, \sigma)$ 
under the assumption that 
$p_{gamb} = p_{true}$.

By a symmetry argument, in the case $p_{gamb} = p_{true} $ we must have 
$[x_1,x_2] = [p_{book} - y, p_{book} + y]$
for some $y$. 
To apply (\ref{Mbook1}),
\[ x_2 - p_{true} = y + \xi; \quad p_{gamb}  + L - x_2 = L - y - \xi .\]
Neglecting odd orders of $\xi$ (which will have expectation zero)
\[ (x_2 - p_{true}) (p_{gamb}  + L - x_2)^2 = y(L-y)^2 + (y - 2(L-y)) \xi^2 \]
and so 
\begin{equation}
 \Ex [ (x_2 - p_{true}) (p_{gamb}  + L - x_2)^2 ] =   y(L-y)^2 + (3y-2L)\sigma^2 .
 \label{Mbook4}
 \end{equation}
Maximizing this over $y$ involves setting $d/dy( \cdot) = 0$ and solving the quadratic equation, yielding 
the optimal value
\[ y^* = 
\sfrac{1}{3} \left(  2 - \sqrt{1 -  \sfrac{9\sigma^2}{L^2} }     \right) L
 \]
 provided $\sigma \le L/3$.
So in this model the bookmaker's strategy is to use the spread interval
\begin{equation}
 [x_1,x_2] = [p_{book} - y^*, p_{book} + y^*] .
 \label{Mbook6}
 \end{equation}
In our case where  $p_{gamb} = p_{true}$,
the value of (\ref{Mbook4}) at $y^*$ works out as
\[
\sfrac{2}{27} \left( 1 + (1- \sfrac{9\sigma^2}{L^2})^{3/2} \right) 
L^3 .
\]
The contribution from the second term in (\ref{Mbook1}) is the same by symmetry, so 
taking account of the pre-factor $\sfrac{\kappa}{4L}$ in (\ref{Mbook1})  we find, 
\begin{equation}
\Ex  [ \mbox{ mean gain to bookmaker   } ] = 
\kappa  h( \sfrac{\sigma^2}{L^2} )  L^2 
\label{Mbook5}
\end{equation}
where
\[ h(u) = 
\sfrac{1}{27} \left( 1 + (1- 9u)^{3/2} \right) .
\]
The function $h(u)$ is shown in Figure \ref{Fig:AB} (left).
The result corresponds to intuition:
 as $\sigma$ increases the gamblers would start to profit  from having the correct consensus probability,
and so  the bookmaker needs to widen the spread interval in response.  
Note also that $h^\prime(0) < 0$, meaning that the cost of the bookmaker's error scales as $\sigma^2$ for small $\sigma$,
 continuing the ``error squared" principle.
 
So the general {\em allowance issue} question
 ``should agents try to adjust their strategy by making allowance for error in estimating probabilities?" 
 has the answer ``yes" for bookmakers in this model, because the optimal spread interval depends on $\sigma$.
 However the practical issue is whether bookmakers could estimate their own $\sigma$ in order to use 
 (\ref{Mbook6}), leading to deeper issues discussed in section \ref{sec:rPT}.

\subsection{Further variations}
\begin{figure}[h]
\caption{Mean gain to bookmaker.  (Left) as a function of bookmaker's error, if gamblers unbiased.
(Right) as a function of gamblers bias, for given bookmaker error.}
\label{Fig:AB}
\includegraphics[width=2.50in]{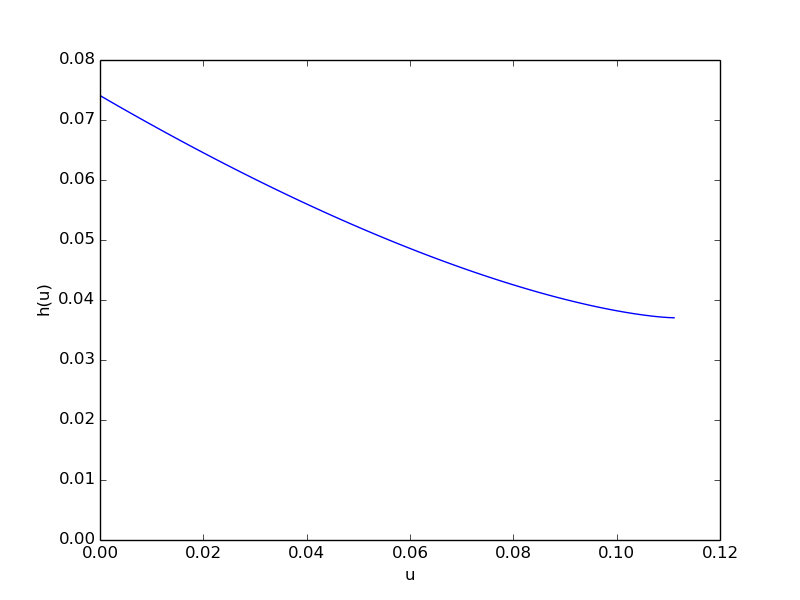}
\includegraphics[width=2.50in]{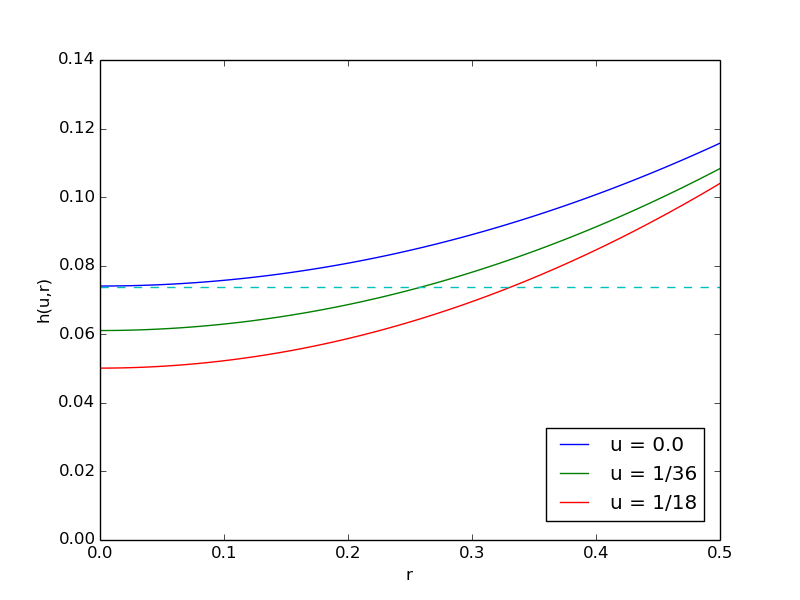}
\end{figure}

The assumption above that $p_{gamb} = p_{true}$ is rather pessimistic from the viewpoint of the bookmaker, who presumably will do better in the more realistic case 
that $\Delta = p_{gamb} - p_{true}$ is non-zero.
How to set the spread interval in that case is one of 
many variants of the models above which may deserve further study.
Here we merely note that
one can calculate the mean gain to bookmaker when the $\Delta = 0$ spread interval (\ref{Mbook6}) 
is used when  $\Delta  \neq 0$; this gives 
%(see section \ref{sec:append}) 
a  formula of the form
  \begin{equation}
\Ex  [ \mbox{ mean gain to bookmaker  } ] = 
\kappa  h^*( \sfrac{\sigma^2}{L^2} , \sfrac{\Delta}{L})  L^2 .
\label{Mbook7a}
\end{equation}
 Figure \ref{Fig:AB}  (right) shows the function $r \to  h^*(u, r)$ for $ u = 0, 1/36, 2/36$.
 As $r$ increases the initial penalty (to bookmaker) of the bookmaker's inaccuracy is offset increasingly by the
 inaccuracy of the gamblers' consensus probability.

\section{Bet I'm better than you!}
\label{sec:skill}
The previous examples involve assessing probabilities of events outside your control. 
A rather different context involves gambling on your own skill.  
We will give two examples.

\subsection{Bradley-Terry games}
Suppose each player has a real-valued ``skill" level $x$, and when player $A$ plays player $B$
\begin{equation}
\Pr(A \mbox{ beats } B) = L(x_A - x_B)
\label{B-T1}
\end{equation}
for a specified function $L$, the conventional choice being the logistic function
\begin{equation}
 L(u) := \frac{e^u}{1 + e^u}, \ - \infty < u < \infty . 
 \label{B-T2}
\end{equation}
This is often called the {\em Bradley-Terry model} and is closely connected with the {\em Elo rating} 
method of estimating skill levels based on past win/lose results -- see \cite{me-Elo} for an introduction.
For our purposes we imagine a game somewhat analogous to a  two-person version of poker, in that games only occur when both players agree to play, and all bets are at even odds -- you win or lose $1$ unit. 
In this setting, if players knew their skills, then bets would never happen 
(in our narrow rationality context\footnote{In reality there are many reasons one might wish to play against a better player, so it may indeed be rational to accept an unfavorable bet to encourage opponent to play.}), because the less skillful player
would refuse to bet..
But by analogy with the simple 2-person bet in section \ref{sec:2pers},
if players don't know the skills exactly then there will be occasions where each player perceives themself as better 
and so are willing to play and bet.

The natural model in our framework is rather complicated, because you may not know your own skill:
so there are 4 errors $\xi_{AA}, \xi_{AB}, \xi_{BA}, \xi_{BB}$ to consider, in the format
\[
\mbox{ A's perception of B's skill} = x_B + \xi_{AB} .
\]
So $A$ is willing to play if $x_A  + \xi_{AA} > x_B + \xi_{AB} $, and $B$ is willing to play if $x_B + \xi_{BB} > x_A + \xi_{BA}$.
So the game occurs if and only if
\[ \xi_{BB}  -\xi_{BA} >   x_A - x_B > \xi_{AB} - \xi_{AA} . \]
Now suppose that as $A$ varies the quantities $ \xi_{AA}  -\xi_{AB} $ are 
independent differently scaled copies of a ``standardized" distribution $\zeta$ 
which is symmetric with variance $1$;
that is
\[  \xi_{AA}  -\xi_{AB} := \sigma_A \zeta_A, \quad  \xi_{BB}  -\xi_{BA} := \sigma_B \zeta_B \]
We now have
\begin{equation}
 \Ex (\mbox{ gain to } A) = \Pr(\sigma_A \zeta_A < x_A-x_B) \Pr(\sigma_B \zeta_B > x_A-x_B)\ (2 L(x_A - x_B) -1). 
 \label{xiAB}
 \end{equation}
 In the context of a given individual $A$ encountering different possible opponents $B$, we need to average over
 some distribution for the skill difference $u = x_A - x_B$, and we will take this as uniform on the real line.\footnote{In Bayesian terms this is an improper prior; because bets only occur when when skill difference is small, we get a ``proper" distribution of actual bets, although (\ref{xiAB2}) should be more  precisely interpreted as a {\em rate} of gain to $A$.}
 So averaging over $B$
\begin{equation}
 \Ex (\mbox{ gain to } A) = \int_{- \infty}^{\infty}  \Pr(\sigma_A \zeta_A < u) \Pr(\sigma_B \zeta_B > u)\ (2 L(u) -1) \ du .
 \label{xiAB2}
 \end{equation}
Here $\sigma_B$ may be random but is assumed independent of skill difference $u$.

\begin{figure}[h!]
\caption{Gain to $A$ as a function of opponent error, for given values of $A$ error.}
\label{Fig:integrate}
\hspace*{1.1in} \includegraphics[width=2.50in]{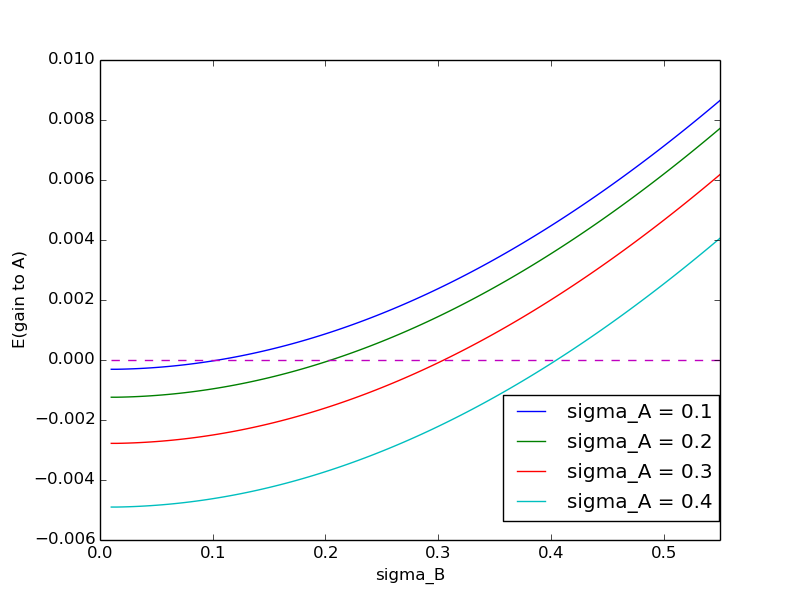}
\end{figure}

We want to determine the effect of the errors in estimating skill -- if the variability $\sigma_A$ of $A$'s estimate is usually smaller than the
variability  $\sigma_B$ of an opponent $B$'s estimate then we anticipate that $A$ will have positive mean gain.  
This is not easy to see from the formula (\ref{xiAB2}).  
If we take $\sigma_B$ to be constant then Figure \ref{Fig:integrate} shows numerical values.
We see the usual error-squared behavior.

\section{Kelly rules}
\label{sec:kelly}

The general {\em Kelly strategy} \cite{kelly},
when  a range of bets (some favorable) with known probability distributions of outcomes are available, is to divide your fortune (normalized to 1 unit) between bets or reserved, and to do so 
in the way such that the random value ($Z$ units) of your fortune after the bets are resolved maximizes 
$\Ex[ \log Z]$.  This value $\Ex [\log Z]$ is the resulting optimal long-term {\em growth rate}.
%Under the Kelly strategy outlined in section \ref{sec:kelly1}, for available bets with known outcome probabilities 
%one can optimize one's portfolio of bets to maximize $\Ex [\log Z]$, which will be the resulting
%optimal {\em growth rate}.\footnote{A portfolio of bets at one step turns initial fortune $1$ into random fortune $Z$.}
But what happens when probabilities are unknown?

In this section we consider the simple setting of betting at even odds, on events with probability close to $0.5$.
If we bet a small proportion $a$ of our fortune and the event occurs with probability 
$0.5 +\delta$ for small $\delta$ then to first order
(we use this approximation throughout)
\begin{equation}
 \mbox{ growth rate } = 2a\delta - a^2/2 .
\label{eq:kelly}
\end{equation}
So for known $\delta > 0$ the optimal choice of proportion is $a = 2 \delta$ 
and the resulting optimal growth rate is $2 \delta^2$. 
Formula (\ref{eq:kelly}) remains true for small $\delta < 0$ but of course here the optimal choice is $a = 0$.

\begin{figure}[h!]
\caption{Growth rate in the Kelly model}
\label{Fig:kelly_1}
\hspace*{0.7in}
\includegraphics[width=3.50in]{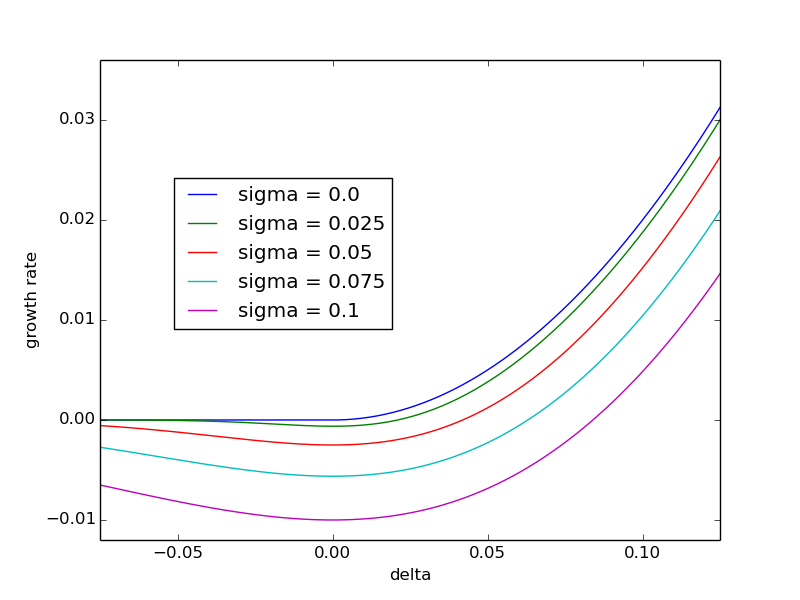}
\end{figure}

In our context there is a perceived probability $0.5 + \delta_{perc}$
and we make the optimal choice based on the perceived probability, that is to bet a proportion
$a = \max(0, 2 \delta_{perc})$. 
We use our usual model for unknown probabilities
\[ \delta_{perc} = \delta_{true} + \xi .
\]
The growth rate $2a\delta_{true} - a^2/2$ can be rewritten as 
\begin{eqnarray*}
 \mbox{ growth rate } &=&  2(\delta_{true}^2 - \xi^2) \ \mbox{ if } \xi > - \delta_{true} \\
 &=& 0 \mbox{ else}.
 \end{eqnarray*}
 Now assume that $\xi$ has Normal($0,\sigma^2$) distribution.
 We can evaluate the expectation of  the growth rate in terms of the pdf $\phi$ and the cdf $\Phi$ of the standard Normal $Z$.
 For $\delta := \delta_{true}$,
 \begin{eqnarray*}
 \Ex [\mbox{ growth rate }] &=&2   \Ex [(\delta^2 - \sigma^2Z^2)1_{(\sigma Z > - \delta)} ]\\
 &=& 2 \left(  \delta^2  \Phi(\delta/\sigma)   - \sigma^2 S(-\delta/\sigma)      \right)
 \end{eqnarray*}
 where
 \[ S(y) := \Ex [Z^2 1_{(Z > y)}] = y \phi(y) + \Phi(-y) . \]
 Putting this together,
 \begin{equation}
 \Ex [\mbox{ growth rate }] = 2 (\delta^2 - \sigma^2)\Phi(\delta/\sigma) + 2 \sigma \delta \phi(\delta/\sigma) .
 \label{Kgrowth}
 \end{equation}
 Figure \ref{Fig:kelly_1} shows the growth rate as a function of $\delta := \delta_{true}$ for several values of $\sigma$.  
 For $\delta > 0$ we see the usual ``quadratic" behavior in both $\delta$ and $\sigma$.
 
This is another setting where the general {\em allowance issue} question
``could agents do better if they knew the typical accuracy of their perceived probabilities and adjusted 
their actions somewhat?" arises.  
In preliminary study we have found it difficult to improve on  the growth rate (\ref{Kgrowth}).

The ``available favorable bets'" requirement for Kelly criterion means it is typically used in a stock market context.
Unlike our previous examples, there is no interaction with another agent, and we classify this as a
``game against nature" as will be discussed in section \ref{sec:nature}.

\section{Pistols at dawn}
\label{sec:pistols}
Let us imagine an 18th century duel: two opponents (A and B) walk toward each other, each with a 
(not very accurate) pistol with a single shot allowed.
As a very simple model, the probability $p(x)$ of a shot incapacitating the opponent depends on the individual and the distance $x$ apart.  If the first shot misses, the opponent can then advance and certainly incapacitate the other.
Assume $A$ is the worse shot, so the hit probabilities are increasing (as distance $x$ decreases) with $p_B(x) \ge p_A(x)$.

In the case of known probability functions, one can readily see that the natural strategy for $A$ is to shoot (if not previously shot at) at distance around $x^*$ 
defined as the solution of 
\[ p_B(x^*) + p_A(x^*) = 1 \]
and also the strategy for $B$ is to shoot (if not previously shot at) at the same distance.\footnote{At {\em approximately} the same distance; we can ignore the possibility of simultaneous shots.}.
So
\begin{equation}
 \Pr(B \ \wins) = p_B(x^*), \ \Pr(A \ \wins) = p_A(x^*) .
 \label{eq:duel1}
\end{equation}
To study the effect of unknown probability functions we need a more detailed model.
Assume that each participant has an accuracy parameter $\rho > 1$ and that $p(x) = \min(\rho/x, 1)$.  
Equation (\ref{eq:duel1}) becomes 
$\rho_A/x + \rho_B/x = 1$ 
and so $x^* = \rho_A + \rho_B$  and so
\begin{equation}
 \Pr(A \ \wins) = \rho_A/(\rho_A + \rho_B); \quad \ \Pr(B\ \wins) =  \rho_B/(\rho_A + \rho_B) .
\label{AAB}
\end{equation}
This is in fact equivalent via transformation\footnote{Under the re-parametrization $\rho \to x = \log \rho$.
But our error model (additive noise $\xi$) is different in the two parametrizations.}
to the Bradley-Terry model (\ref{B-T1}, \ref{B-T2}),
but in this context participants are compelled to participate rather than having the choice to bet or not to bet.
Consider the model where participants know their own ability but not the opponent's. 
That is,

$A$ perceives B's  accuracy parameter as $\rho_B + \xi_A$; 

and  $B$ perceives A's accuracy parameter as $\rho_A + \xi_B$.

\noindent
So $A$ plans to shoot at distance  $x_A^* = \rho_A + \rho_B + \xi_A$;
and $B$  plans to shoot at  distance $x_B^* = \rho_A + \rho_B + \xi_B$.
The outcome of the duel is
\begin{itemize}
\item if $\xi_A > \xi_B$ then $A$ shoots first and $A$ wins with probability 
$\rho_A/(\rho_A + \rho_B + \xi_A)$;
\item  if $\xi_A < \xi_B$ then $B$ shoots first and so $A$ wins with probability 
$1 - \rho_B/(\rho_A + \rho_B + \xi_B)$.
\end{itemize}
Because of the rather arbitrary parametrization, in this model there is a first order effect -- 
if $\xi_A$ and $\xi_B$ have mean zero and the same small variance $\sigma^2$, then
the change in probability that $A$ wins scales generically as $\sigma$ rather than $\sigma^2$.
For a  more robust observation, if $A$ knows $B$'s ability but has the opportunity to mis-represent $B$'s knowledge of $A$'s ability, which direction of mis-representation is desirable?
In this model $A$ would like $B$  to over-estimate $A$'s ability ($\xi_B > 0$), which would motivate $B$ to shoot earlier and less accurately.

\section{How valuable is it really?}
\label{sec:Gumbel}
The models discussed so far have featured discrete outcomes with unknown probabilities.
More generally we can consider continuous outcomes, that is a range of possible utilities, where one's perception of the utility is inaccurate. 
Here we study one specific mathematical model chosen to be analytically tractable,
and to which two slightly different stories can be attached. 
Note that in this kind of setting, the details of the model can make a large difference.

Perhaps the simplest type of such model is the ``choose the best item from many items" type.

\subsection{Choosing the best item}
Here one agent (you) needs to choose one out of a given set of items.
In our model,  the true utilities of the items are distributed as the inhomogeneous  Poisson process of intensity function 
\begin{equation}
 \lambda(x) := e^{-x},  \ - \infty < x < \infty 
 \label{def:lambda}
 \end{equation}
 this distribution arising from classical extreme value theory \cite{resnick}.
 In that model the largest true utility  $X_{(1)}$ has the Gumbel distribution with c.d.f. and p.d.f.
\[ G(x) = \exp(- e^{-x}),
\ g(x) = e^{-x} \exp(- e^{-x}),
 \ -\infty < x < \infty . \]
If you knew the true utilities then you would pick the item and gain utility $X_{(1)}$.
However, in our model your perceived utility of an item of true utility $x$ is $x + \xi$ for i.i.d. random $\xi$.
So by choosing the item with largest perceived utility, you gain some utility $Y$ which may or may not 
be  $X_{(1)}$, and so you incur a ``cost" $X_{(1)} - Y \ge 0$. 
We will study the mean cost.

Note first that the mean difference between the largest two true utilities can be 
calculated in terms of 
\[ N(x) := \mbox{ number of items with true utility }>x \]
via
\begin{eqnarray}
\Ex  [ X_{(1)} - X_{(2)}  ] &=& \int_{-\infty}^\infty \Pr(N(x) = 1) \ dx \label{Nx1} \\
&=&  \int_{-\infty}^\infty e^{-x}  \exp(- e^{-x}) \ dx \nonumber \\
&=& 1 . \nonumber
\end{eqnarray}
Note also that replacing $\xi$ by $\xi - \Ex[\xi]$ makes no difference, so we may assume $\Ex [\xi] = 0$.

Now suppose the errors $\xi_\sigma$ have Normal($0,\sigma^2$) distribution and write
$\phi_\sigma$ and $\Phi_\sigma$ for the p.d.f and c.d.f. of that distribution.
The process of pairs 
(true utility of item, perceived utility of item) 
is Poisson with intensity 
\[
\lambda_2(x,y) = e^{-x} \phi_\sigma(y-x), \ - \infty < x, y < \infty .
\]
We can calculate the mean cost in the style of (\ref{Nx1}) as follows.
\[
\Ex [ \mbox{ cost } ] =
\int_{-\infty}^\infty  \int_{-\infty}^{x_{(1)}} \int_{-\infty}^{\infty}
g(x_{(1)}) \Phi_\sigma(y-x_{(1)})      \ \Pr(M_{(x,x_{(1)})} < y) \ \Pr( M_{(-\infty, x]} \in dy) \ dx dx_{(1)} 
\]
where $M_I$ is the maximum perceived utility amongst items with true utility in interval $I$.
And
\begin{eqnarray*}
\Pr(M_I < y) &=& 
\exp \left( - \int_I e^{-u}(1 - \Phi_\sigma(y-u)) \ du \right) \\
\Pr(M_I \in dy) &=& \Pr(M_I < y) \ \cdot \ \int_I e^{-u} \phi_\sigma (y-u) \ du. 
\end{eqnarray*}

\begin{figure}[h!]
\caption{Mean cost in the ``Choosing the best item" model.}
\label{Fig:exp_exp}
\hspace*{0.6in}
\includegraphics[width=3.50in]{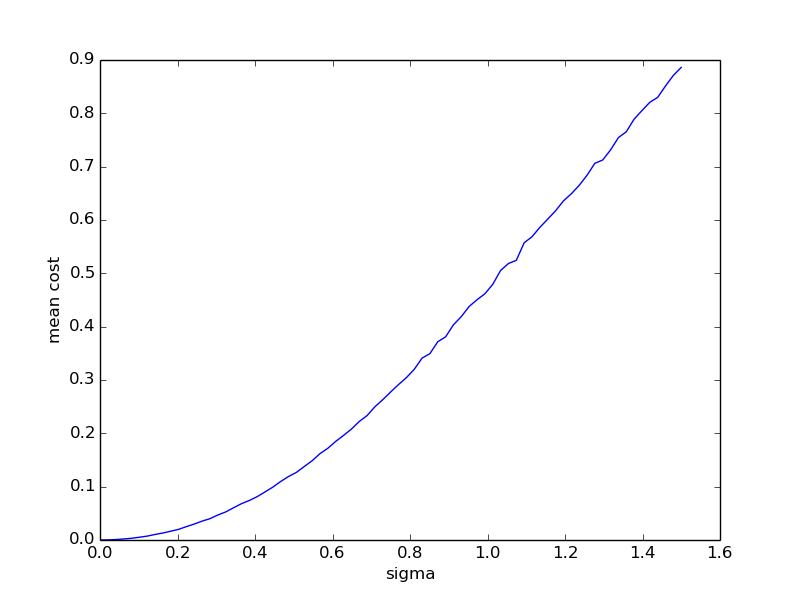}
\end{figure}

Figure \ref{Fig:exp_exp} shows numerical values. 
%{\tt xxx This looks analytically intractable for Normal $\xi$ but maybe can get something for Laplace distribution ???}
Note this is what we call a  ``game against nature" and  
(as will be discussed in section \ref{sec:nature}) what we find is the expected order
$\sigma^2$  cost for $\sigma \ll 1$.
The subsequent linear growth reflects the model assumption of exponentially growing numbers of items of decreasing utility.

\subsection{Profit or loss at auction?}
We can re-use the mathematical model at (\ref{def:lambda}) as an auction model, as follows.  
Now there is one item being sold, and a large number of agents  who can bid.
For each agent, the future benefit of acquiring the item is some value $x$ which the agent does not know exactly,
but instead perceives the value as $y = x + \xi$ for random $\xi$, i.i.d. over agents, and bids that perceived value.
In this model the points of the Poisson process (\ref{def:lambda})  represent the actual values $x_i$ of the given item to the different agents.
The winner of the auction is the agent $i$ whose bid $y_i = x_i + \xi_i$ is largest.
In a traditional sealed-bid auction, that winner pays their own bid amount $y_i$ and so makes a profit
of $ x_i - (x_i + \xi_i) =  - \xi_i$.  
In a modern {\em Vickrey auction} the winner pays the second-highest bid amount $y_{(2)}$, to mimic a live auction, 
and so makes a profit of $ x_i  - y_{(2)}$.

\begin{figure}[h!]
\caption{Mean gain in the auction model.}
\label{Fig:auction}
\hspace*{0.6in}
\includegraphics[width=3.50in]{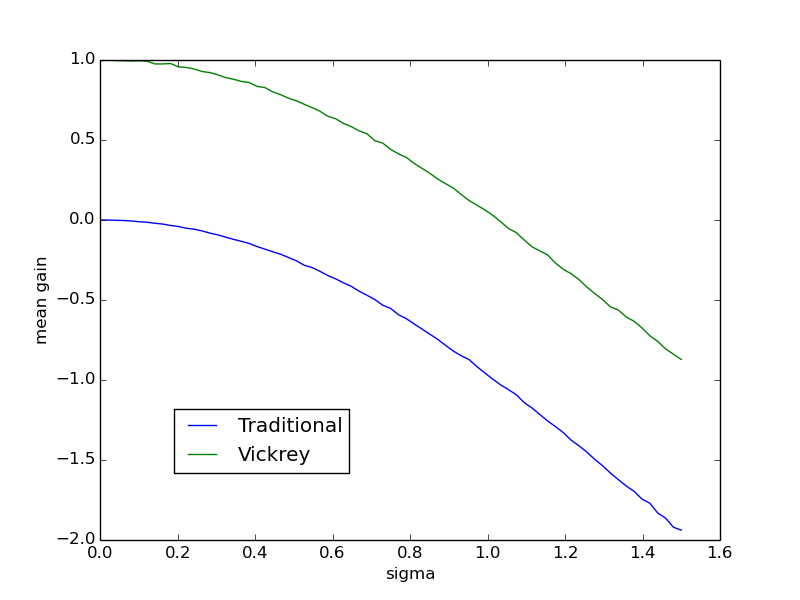}
\end{figure}

Figure \ref{Fig:auction} shows  numerical values of the mean gain as a function of $\sigma$, when $\xi$ has Normal($0,\sigma^2$) distribution. 
As usual we see $\sigma^2$ behavior for small $\sigma$.
It is curious that the difference in mean gains between the two auction protocols remains roughly constant as $\sigma$ varies.

{\bf Game theory aspects.}
In the auction model above, agents simply bid their perceived value without attention to competing bidders.
There has been extensive theoretical and applied work on game-theoretic aspects of auctions, in particular in the context 
of Google auctions for sponsored advertisements  \cite{sponsored}.
Chapter 14 of \cite{karlin} gives an introduction to relevant theory.
The reader might consider how to combine our framework with game-theoretic models.

\section{Background and discussion}
\label{sec:discuss}

\subsection{Brief review of prediction tournaments}
\label{sec:rPT}
This is the closest topic that has been well studied by elementary mathematics, and so is natural background for this article. 
A {\em prediction tournament} (see \cite{mePT} for an introductory account) consists of a collection of questions of the form 
``state a probability for a specified real-world event
happening before a specified date".\footnote{In actual tournaments one can update probabilities as time passes, 
but for simplicity we  consider only a single probability prediction for each  question, and only binary outcomes.}
Typical questions in a current (August 2020) tournament are: before 1 January 2021,
\begin{itemize}
\item  will Benjamin Netanyahu cease to be the prime minister of Israel ?
\item  will there be a lethal confrontation between the national military or law enforcement forces of Iran and Saudi Arabia?
\item  will the Council of the European Union request the consent of the European Parliament to conclude a European Union-United Kingdom trade agreement?
\end{itemize}
These are unique events, not readily analyzed algorithmically.  
In the popular book \cite{fox}
contestants are advised to combine what information they can find about the specific event 
with some ``baseline" frequency of roughly analogous previous events.
Scoring is by squared error: if you state probability $q$ then on that question then
\[ \mbox{
score = $(1-q)^2$ if event happens;  score = $q^2$ if not. 
}\]
Your tournament score is the sum of scores on each question.
As in golf one seeks a {\em low} score. 
Also as in golf, in a {\em tournament} all contestants address the same questions; 
it is not a single-elimination tournament as in tennis.

\paragraph{A key starting insight.}
The starting insight (known for a long time but apparently not widely known) 
is that in a prediction tournament
\begin{quote}
even though true probabilities are completely unknown, one can determine 
(up to small-sample chance variability) the  {\em relative} abilities of contestants 
at estimating true probabilities.
\end{quote}
That is, with unknown true probabilities $(p_i)$,  if you announce probabilities $(q_i)$
then by elementary algebra \cite{mePT} the true expectation of your score equals
\[ \sum_i p_i(1-p_i) + \sum_i (q_i - p_i)^2 . \]
The first term is the same for all contestants, so if $S$ and $\hat{S}$ are the tournament scores for you 
and another contestant in an $n$-question tournament, then
\begin{equation}
n^{-1/2} (\Ex S - \Ex \hat{S}) = \sigma^2 - \hat{\sigma}^2 
\label{PT1}
\end{equation}
where
\[ \sigma : = \sqrt{n^{-1} \sum_i (q_i - p_i)^2 } \]
is your RMS error in predicting probabilities and $\hat{\sigma}$ is the  other contestant's RMS error.
Thus by looking at differences in scores one can, in the long run, use (\ref{PT1}) to estimate {\em relative} abilities at prediction,
as measured by RMS error of predicted probabilities.

The ``insight" here is the clear relevance of RMS error (rather than some other measure of size of error) in this context.

Of course ``long run" assertions deserve some consideration of non-asymptotics. 
Under fairly plausible more specific assumptions, \cite{mePT} shows that in a 100-question tournament,
if you are $5\%$ more accurate than me (e.g. your RMS error is 15\% while mine is 20\%),
then your chance of beating me is around 75\%; increasing to around 92\% if a 10\% difference in RMS error.
So there is quite a lot of chance variability due to event outcomes.

Whether one can estimate $\sigma$ itself, that is a contestant's {\em absolute} rather than {\em relative} error, is a deep question.
Under a certain model, it is shown in \cite{mePT} that from the distribution of scores in a tournament one can determine roughly the $\sigma$ for the best-scoring contestants, but the assumptions of the model
are not readily verifiable.

\subsection{On unbiased and calibrated estimates}
\label{sec:calibrate}
We have not attempted  to model how an agent's perceived probabilities are obtained -- just some ``black box" method. 
However an agent can record predictions and outcomes to check whether the {\em calibrated} property \cite{tetlock}
\begin{quote}
 the long-run proportion of events with perceived probability near $p$ is  near $p$
 \end{quote}
 holds.   In principle one can mimic this property in the long run: if one finds that only 15\% of events 
 for which one predicted 20\% actually occur, then in future when one's ``black box" outputs 20\%,
 one instead predicts $15\%$.
 So {\em calibration} is not so unrealistic for serious gamblers who monitor their performance.  
 
 Now imagine a scatter diagram of $(p_{perc}, p_{true})$. 
 Freshman statistics reminds us that there are two different regression lines
 (for predicting one variable given the other).
 The calibrated property is that one of those lines is the diagonal on $[0,1]^2$;
 our {\em unbiased} property $\Ex \xi = 0$ is that the other line is the diagonal.
 These are logically different, but
 if the errors $\xi$ are small then the lines are not very different\footnote{They differ more near $0$ or $1$, but recall we are implicitly dealing with probabilities not near $0$ or $1$.}; 
 so {\em calibrated} and {\em unbiased}  
capture the same intuitive idea. 
 
Note also that  for simplicity we have usually assumed that the variance of $p_{perc}$ for given $p_{true}$ is a constant $\sigma^2$ 
 not depending on $p_{true}$, but without that assumption results such as (\ref{eq:2-person}, \ref{bet-PM})  remain true by replacing $\sigma^2$ by the conditioned variance.

\subsection{The error-squared principle in bets against nature}
\label{sec:nature}
We have emphasized bets against humans because, in a sense described below, ``bets agains nature" 
are simpler.
One type of 
``decisions under unknown probabilities" setting can be abstracted as follows.
You have an action which involves choosing a real number $x$.
The outcome depends on whether an event of unknown probability $p$ occurs. 
The mean gain (in utility) is a function
$\mathbf{gain}(x,p)$.
If you knew the true probability $p_{true}$ then you could choose the optimal action 
\[ x_{opt} = \arg \max_x \mathbf{gain}(x,p_{true}) . \]
Instead you have a perceived probability $p_{perc}$ and so you choose an actual action
\[ x_{actual} = \arg \max_x \mathbf{gain}(x,p_{perc}) . \]
So the cost of your error is the difference
\[ \mathbf{cost}(p_{perc}) := \mathbf{gain}(x_{opt},p_{true}) - \mathbf{gain}(x_{actual},p_{true}) . \]
And then  calculus tells us 
 that for a generic smooth function $\mathbf{gain}$
 we have $|x_{actual} - x_{opt}|$ is order $|p_{perc} - p_{true}|$ and so 
 because $x \to \mathbf{gain}(x,p_{true}) $ is maximized at $x_{opt}$,
\[ \mathbf{cost}(p_{perc}) \mbox{ is order } (p_{perc} - p_{true})^2 \mbox{ as }  p_{perc} \to p_{true} . \]
So the {\em cost} of an error in estimating a probability is typicaly scales as the square of the {\em size} of the error.

Let's see how this arises in more detail in one simple example.

{\bf Example.}
Suppose you are planning a wedding, several months ahead, and need to choose now between
an outdoor venue A, which you would prefer if you knew it would not be raining, 
or an indoor venue B, which you would prefer if you knew it would be raining.
How to choose?  
Invoking utility theory (and ignoring the difficulty of actually assigning utilities), there are 4 utilities
\begin{itemize}
\item (choose A): utility $=a$ if no rain, utility $=b$  if rain
\item (choose B): utility $=c$ if no rain, utility $=d$  if rain
\end{itemize}
where $a>c$ and $d>b$.  
Now  you calculate the expectation of the utility in terms of the probability $p_{true}$ of rain:
\begin{itemize}
\item  (choose A): expected utility $ = p_{true}b + (1-p_{true})a$
\item  (choose B): expected utility $ = p_{true}d + (1-p_{true})c $
\end{itemize}
There is a critical value $ p_{crit}$  where these mean payoffs are equal, and this is the solution of 
$\frac{p_{crit}}{1-p_{crit}} = \frac{a-c}{d-b}$.
Note that $p_{crit} $ depends only on the utilities.
If you knew $p_{true}$ your best strategy would be
\[ \mbox{
do action A if $p_{true} < p_{crit} $, do action B if $p_{true} > p_{crit} $}.
\]
Instead all we have is our guess $p_{perc}$, so we use this strategy 
but based on $p_{perc}$ instead of $p_{true}$.

What is the cost of not knowing $p_{true}$?
If $p_{true}$  and $p_{perc}$ are on the same side of $p_{crit} $ then you take the 
optimal action and there is zero cost; 
if they are on opposite sides then you take the sub-optimal action and the cost is 
\[
\mbox{ $|p_{true} - p_{crit} |  z$ 
where $z = a- b - c + d > 0$.} 
\]
Note that $z$ and $p_{crit}$ are determined by the parameter $\theta := (a,b,c,d)$ which we take to be of order $1$.
 In our setting, 
 $(\theta,p_{true},p_{perc})$ come from some smooth distribution on triples.
 As in the previous example, 
  the event that $ \cost > 0$, and then the  mean of $ \cost $ given $ \cost  > 0$,
 are both order $|p_{perc} - p_{true}|$, and so 
 \[ \Ex  [\cost ] \mbox{ is order } (p_{perc} - p_{true})^2 . \]

\section{Remarks on the academic literature}
\label{sec:sophist}
We do not know any comparable elementary discussion of our general topic -- 
-- quantitative study of consequences of the fact that numerical probabilities can often only be roughly estimated -- but the topic
has certainly been considered at a more sophisticated level in many contexts, some of which we describe here.

Of course  one can estimate an unknown probability
in the classical context of  ``repeatable experiments".
In the standard mathematical treatment of stochastic processes models (e.g. \cite{pinsky}) 
one has a model with several  parameters (some relating to probabilities) 
and the implicit relevance to the real-world is that, with enough data, one can estimate parameters
and thereby make predictions.
The recent field of 
{\em uncertainty quantification} \cite{uncertain} seeks to address all aspects of error within complicated models,
but implicitly assuming the model is known to be scientifically accurate or can be calibrated via data.
Our focus, exemplified by the prediction tournament examples from section \ref{sec:rPT} or sports gambling, 
is on contexts of ``unique events" where one has no causal models and directly relevant  past data, 
necessitating human judgment not algorithms .

Sports gambling is too huge a field to survey here; see \cite{sports_betting} for a non-mathematical introduction, and  browse  {\em Journal of Quantitative Analysis in Sports} for articles implicitly relating to events one could bet upon.

Already mentioned is  the detailed account of actual bookmaker strategy in the U.S. sports gambling context
by Levitt  \cite{levitt}.
For one aspect of gambling on horse races to illustrate the style of literature,  Green et al.  \cite{midas} seek to explain the observed 
low return from betting on longshots in parimutuel markets as follows.
\begin{quote}
The track deceives naive bettors by suggesting inflated probabilities for longshots and depressed probabilities for favorites. This deception induces naive bettors to underbet favorites, which creates arbitrage opportunities for sophisticated bettors and, from that arbitrage, incremental tax revenue for the track.
\end{quote}

There has been considerable theoretical work on the
{\em forecast aggregation problem},  that is how to combine different probability forecasts 
 into a single consensus forecast.
See \cite{ariel,cross,pemantle} for representative recent work.

\subsection{Errors in estimating near-critical parameters}
\begin{figure}[h!]
\caption{Pandemic model analysis from \cite{Castro26190} .}
\label{Fig:pandemic}
\hspace*{0.6in}
\includegraphics{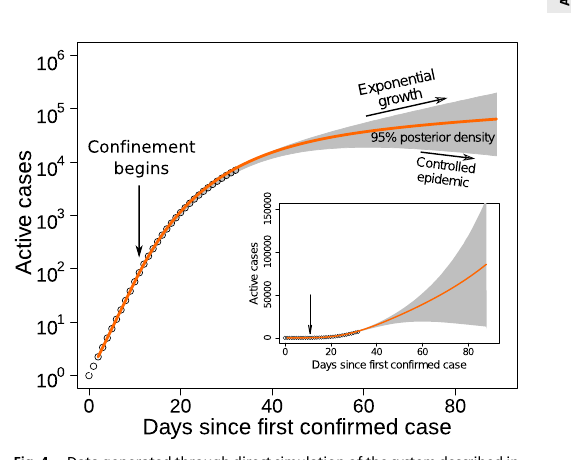}
\end{figure}

One setting where our squared-error paradigm is (rather obviously) not applicable 
concerns models with a parameter which has a critical value 
(meaning qualitatively different behavior on the two sides of the parameter).
If the true value is near the critical value then uncertainty about the value may directly lead 
to great uncertainty about the behavior.
The most elementary textbook case  concerns sub- or super-criticality  in branching processes.
The 2020 COVID-19 pandemic led to extensive analysis of epidemic models. 
Figure \ref{Fig:pandemic} is a graphic from \cite{Castro26190} who emphasize that, within a 
fairly plausible Bayes model, whether or not a confinement policy will control the epidemic cannot be 
determined even with complete knowledge up until confinement begins.

\end{document}